\documentclass[12pt]{article}
\usepackage{amsmath,amssymb, amsfonts, amsthm,amscd}
\usepackage[T2A]{fontenc}
\usepackage[cp1251]{inputenc}
\usepackage[russian,ukrainian,english]{babel}
\usepackage{geometry}

\newcommand{\case}[1]{\textbf{Case $\mathbf{#1}$.}\ }

\newtheorem{Theorem}{Theorem}
\newtheorem{Corollary}{Corollary}

\newtheorem{Proposition}{Proposition}

\newenvironment{TheoremProof}{\textbf{Proof. }}{\par\noindent\textbf{The Theorem is proved.}}

\title{{\Large \textbf{An inequality for the number of vertices with an interval spectrum in edge labelings of regular graphs}}}

\author{\normalsize N.N. Davtyan$^1$, R.R. Kamalian$^2$}

\date{}

\begin{document}

\maketitle

$\\^1$Ijevan Branch of Yerevan State University, e-mail:
nndavtyan@gmail.com $\\^2$The Institute for Informatics and
Automation Problems of NAS RA, \\e-mail: rrkamalian@yahoo.com
\bigskip \bigskip

\begin{abstract}
We consider undirected simple finite graphs. The sets of vertices
and edges of a graph $G$ are denoted by $V(G)$ and $E(G)$,
respectively. For a graph $G$, we denote by $\delta(G)$ and
$\eta(G)$ the least degree of a vertex of $G$ and the number of
connected components of $G$, respectively. For a graph $G$ and an
arbitrary subset $V_0\subseteq V(G)$ $G[V_0]$ denotes the subgraph
of the graph $G$ induced by the subset $V_0$ of its vertices. An
arbitrary nonempty finite subset of consecutive integers is called
an interval. A function $\varphi:E(G)\rightarrow
\{1,2,\dots,|E(G)|\}$ is called an edge labeling of the graph $G$,
if for arbitrary different edges $e'\in E(G)$ and $e''\in E(G)$, the
inequality $\varphi(e')\neq \varphi(e'')$ holds. If $G$ is a graph,
$x$ is its arbitrary vertex, and $\varphi$ is its arbitrary edge
labeling, then the set $S_G(x,\varphi)\equiv\{\varphi(e)/ e\in E(G),
e \textrm{ is incident with } x$\} is called a spectrum of the
vertex $x$ of the graph $G$ at its edge labeling $\varphi$. If $G$
is a graph and $\varphi$ is its arbitrary edge labeling, then
$V_{int}(G,\varphi)\equiv\{x\in V(G)/\;S_G(x,\varphi)\textrm{ is an
interval}\}$. For an arbitrary $r$-regular graph $G$ with $r\geq2$
and its arbitrary edge labeling $\varphi$, the inequality
$$
|V_{int}(G,\varphi)|\leq\bigg\lfloor\frac{3\cdot|V(G)|-2\cdot\eta(G[V_{int}(G,\varphi)])}{4}\bigg\rfloor.
$$
is proved.
\bigskip


Math. Classification: 05C15, 05C78
\end{abstract}

We consider undirected simple finite graphs. The sets of vertices
and edges of a graph $G$ are denoted by $V(G)$ and $E(G)$,
respectively. For a graph $G$, we denote by $\delta(G)$ the least
degree of a vertex of $G$. For any graph $G$ we define a parameter
$\eta(G)$ by the following way: if $G$ is empty then
$\eta(G)\equiv0$, otherwise $\eta(G)$ is equal to the number of
connected components of $G$. If $G$ is a graph, $x\in V(G)$, $y\in
V(G)$, then $\rho_G(x,y)$ denotes the distance between the vertices
$x$ and $y$ in $G$. If $G$ is a graph, $x\in V(G)$, and
$V_0\subseteq V(G)$, then $\rho_G(x,V_0)$ denotes the distance in
the graph $G$ between its vertex $x$ and the subset $V_0$ of its
vertices. For a graph $G$ and an arbitrary subset $V_0\subseteq
V(G)$ $G[V_0]$ denotes the subgraph of the graph $G$ induced by the
subset $V_0$ of its vertices.

For any graph $G$ and its arbitrary subgraph $H$, let us define the
subgraph $Surr[HinG]$ of the graph $G$ as follows:
\begin{displaymath}
\begin{array}{l}
V(Surr[HinG])\equiv\{x\in V(G)/\;\rho_G(x,V(H))\leq1\},\\
E(Surr[HinG])\equiv E(H)\cup\{(x,y)\in E(G)/\;x\in
V(Surr[HinG])\backslash V(H), y\in V(H)\}.
\end{array}
\end{displaymath}

An arbitrary nonempty finite subset of consecutive integers is
called an interval. A function $\varphi:E(G)\rightarrow
\{1,2,\dots,|E(G)|\}$ is called an edge labeling of the graph $G$,
if for arbitrary different edges $e'\in E(G)$ and $e''\in E(G)$, the
inequality $\varphi(e')\neq \varphi(e'')$ holds. For a graph $G$,
the set of all its edge labelings is denoted by $\tau(G)$.

If $G$ is a graph, $x\in V(G)$, $\varphi\in\tau(G)$, then the set
$S_G(x,\varphi)\equiv\{\varphi(e)/ e\in E(G), e \textrm{ is incident
with } x$\} is called a spectrum of the vertex $x$ of the graph $G$
at its edge labeling $\varphi$. If $G$ is a graph,
$\varphi\in\tau(G)$, then $V_{int}(G,\varphi)\equiv\{x\in
V(G)/\;S_G(x,\varphi)\textrm{ is an interval}\}$. The terms and
concepts which are not defined can be found in \cite{West1}.

An upper bound for the cardinality of the set $V_{int}(G,\varphi)$
is obtained in that cases when $G$ is a regular graph with
$\delta(G)\geq2$ and $\varphi\in\tau(G)$.

First we recall the following

\begin{Proposition}\cite{Luhansk17}\label{Thm1}
Let $G$ be a graph with $\delta(G)\geq 2$. Let $\varphi\in\tau(G)$
and $V_{int}(G,\varphi)\neq\emptyset$. Then $G[V_{int}(G,\varphi)]$
is a forest, each connected component of which is a simple path.
\end{Proposition}

\begin{Theorem}\label{Thm2}
If $G$ is a $r$-regular graph, $r\geq2$, $\varphi\in\tau(G)$, then
$$
|V_{int}(G,\varphi)|\leq\bigg\lfloor\frac{r\cdot|V(G)|-2\cdot\eta(G[V_{int}(G,\varphi)])}{2\cdot(r-1)}\bigg\rfloor.
$$
\end{Theorem}

\begin{TheoremProof}
Let $\eta(G[V_{int}(G,\varphi)])=k$.

\case{1} $V_{int}(G,\varphi)=\emptyset$.

In this case the required inequality is the following evident one:
$$
0\leq\bigg\lfloor\frac{r\cdot|V(G)|}{2\cdot(r-1)}\bigg\rfloor.
$$

\case{2} $V_{int}(G,\varphi)\neq\emptyset$.

In this case $k\geq1$. Since $\delta(G)=r\geq2$, then, by the
proposition \ref{Thm1}, $G[V_{int}(G,\varphi)]$ is a forest with $k$
connected components, each of which is a  simple path.

Let $P_1,\dots,P_k$ be all connected components of the forest
$G[V_{int}(G,\varphi)]$.

It is not difficult to see that for $\forall i$, $1\leq i\leq k$,
the equality $|E(Surr[P_iinG])|=(r-1)\cdot|V(P_i)|+1$ holds.

Let us also note that (if $k\geq2$) for arbitrary integers $i'$ and
$i''$ satisfying the inequality $1\leq i'<i''\leq k$, the relation
$E(Surr[P_{i'}inG])\cap E(Surr[P_{i''}inG])=\emptyset$ holds.

Taking into account the evident relation $(\bigcup_{i=1}^k
E(Surr[P_iinG]))\subseteq E(G)$, we obtain
$$
|E(G)|=\frac{r\cdot|V(G)|}{2}\geq\bigg|\bigcup_{i=1}^k
E(Surr[P_iinG])\bigg|=\sum_{i=1}^k |E(Surr[P_iinG])|=
$$
$$
=\sum_{i=1}^k
((r-1)\cdot|V(P_i)|+1)=k+(r-1)\cdot\sum_{i=1}^k|V(P_i)|=k+(r-1)\cdot|V_{int}(G,\varphi)|,
$$
$$
|V_{int}(G,\varphi)|\leq\frac{1}{r-1}\cdot\bigg(\frac{r\cdot|V(G)|}{2}-k\bigg)=\frac{r\cdot|V(G)|-2k}{2\cdot(r-1)}.
$$

Consequently,
$$
|V_{int}(G,\varphi)|\leq\bigg\lfloor\frac{r\cdot|V(G)|-2k}{2\cdot(r-1)}\bigg\rfloor.
$$
\end{TheoremProof}

\begin{Corollary}
If $G$ is a $r$-regular graph, $r\geq2$, $\varphi\in\tau(G)$, then
$$
|V_{int}(G,\varphi)|\leq\bigg\lfloor\frac{r\cdot|V(G)|-2}{2\cdot(r-1)}\bigg\rfloor.
$$.
\end{Corollary}

\begin{Corollary}
If $G$ is a cubic graph, $\varphi\in\tau(G)$, then
$$
|V_{int}(G,\varphi)|\leq\bigg\lfloor\frac{3\cdot|V(G)|-2\cdot\eta(G[V_{int}(G,\varphi)])}{4}\bigg\rfloor.
$$.
\end{Corollary}

\begin{Corollary}
If $G$ is a cubic graph, $\varphi\in\tau(G)$, then
$$
|V_{int}(G,\varphi)|\leq\bigg\lfloor\frac{3\cdot|V(G)|-2}{4}\bigg\rfloor.
$$.
\end{Corollary}

\end{document}